\documentclass[12pt]{amsart}
\usepackage{amsmath, amsthm, amssymb, amscd, mathrsfs, eucal, epsfig}
\usepackage{algorithm} 
\usepackage{fullpage}
\usepackage[noend]{algpseudocode}
\usepackage{enumitem}
\usepackage{pinlabel}
\usepackage{tikz-qtree}
\usepackage{anyfontsize}

\newtheorem{newthm}{Theorem}
\newtheorem{theorem}{Theorem}[section]
\newtheorem{lemma}[theorem]{Lemma}
\newtheorem{proposition}[theorem]{Proposition}
\newtheorem{corollary}[theorem]{Corollary}

\newtheorem*{question}{\bf Question}

\newtheorem{remark}[theorem]{Remark}
\def\beginp{ { \em Proof.} }

\numberwithin{equation}{section}

\def\mystrut{{\rule[-2ex]{0ex}{4.5ex}{}}}

\usepackage{enumitem}

\usepackage{tikz}

\newcommand{\REFEQN}[1] { \begin{equation}\label{#1} }
\newcommand{\ENDEQN}{\end{equation}}
\newcommand{\REFTHM}[1] { \begin{theorem}\label{#1} }
\newcommand{\ENDTHM}{\end{theorem}}
\newcommand{\REFNTH}[1] { \begin{newthm}\label{#1} }
\newcommand{\ENDNTH}{\end{newthm}}
\newcommand{\REFPROP}[1]{\begin{proposition}\label{#1} }
\newcommand{\ENDPROP}{\end{proposition} }
\newcommand{\REFLEM}[1]{\begin{lemma}\label{#1} }
\newcommand{\ENDLEM}{\end{lemma} }
\newcommand{\REFCOR}[1]{\begin{corollary}\label{#1} }
\newcommand{\ENDCOR}{\end{corollary} }

\def\g{\gamma}
\def\G{\Gamma}

\def\ep{\varepsilon}

\def\La{\Lambda}

\def\R{\mbox{$\mathbb R$}}

\def\Z{\mbox{$\mathbb Z$}}

\def\qed{\hfill{q.e.d.}}

\begin{document}

\title{On balanced planar graphs, following W. Thurston}
\author{Sarah Koch}
\address{Department of Mathematics \\ University of Michigan \\ 530 Church Street\\
Ann Arbor MI 48109\\USA}
\email{kochsc@umich.edu}
\author{Tan Lei}
\address{Facult\'e des sciences\\ LAREMA\\ Universit\'e d'Angers\\ 2 Boulevard Lavoisier \\ 49045 Angers cedex \\
France}
\email{tanlei@math.univ-angers.fr}
\date{\today}

%
%
%
%

 \begin{abstract} 
Let $f:S^2\to S^2$ be an orientation-preserving branched covering map of degree $d\geq 2$, and let $\Sigma$ be an oriented Jordan curve passing through the critical values of $f$. Then $\Gamma:=f^{-1}(\Sigma)$ is an oriented graph on the sphere.  In a group email discussion in Fall 2010, W. Thurston introduced {\em balanced planar graphs}  and showed that they combinatorially characterize all such $\Gamma$, where $f$ has $2d-2$ distinct critical values. We give a detailed account of this discussion, along with some examples and an appendix about Hurwitz numbers.
  \end{abstract}
  
  \maketitle

\section{Introduction and Results} 

  Let $S^2$ denote the topological 2-sphere, equipped with an orientation. Let $f:S^2\to S^2$ be an orientation-preserving branched covering map of degree $d\geq 2$, let $R_f$ be the set of ramification values (or critical values) of $f$, and let $C_f$ be the set of critical points of $f$. We say that $f$ is {\em generic} if $|R_f|=2d-2$; note that if $f$ is a generic branched cover, then $|C_f|=2d-2$ as well. Furthermore, every critical point of a generic map $f$ must be simple. 

  Let $\Sigma\subseteq S^2 $ be an oriented Jordan curve running through the $n:=2d-2$ critical values of $f$, where we consider $\Sigma$ to be an oriented graph with vertex set $V_\Sigma=R_f$; each vertex has valence $2$. Then the set $f^{-1}(\Sigma)\subseteq S^2$ is an oriented graph in the domain, with vertex set $V_{f^{-1}(\Sigma)}=f^{-1}(R_f)$.  Among the $nd-n$ elements of $V_{f^{-1}(\Sigma)}$ are: 

\begin{itemize}
\item the vertices with valence $4$, which comprise the set of critical points $C_f$, and 
\item the vertices with valence $2$, which comprise the set $f^{-1}(R_f)-C_f$, often referred to as {\em cocritical points}. 
\end{itemize}
  Forgetting the valence $2$ vertices in $f^{-1}(\Sigma)$  gives an oriented graph $\G$ with vertex set $V_\Gamma=C_f$. We call $\G$ the {\em underlying $4$-valent graph} of $f^{-1}(\Sigma)$. We will refer to a connected component of $S^2-\G$ as a {\em face} of $\G$. 
 

  The question we address in this article is: which oriented $4$-valent  graphs on $S^2$ can be realized as the underlying graph of
$f^{-1}(\Sigma)$ for some  pair $f$ and $\Sigma$? This was asked by W. Thurston as  part of a group email discussion, which he initiated in Fall 2010, and which was motivated by his visionary  question:
\begin{question}[W. Thurston, 2010]
What is the shape of a rational map?
\end{question}
  Throughout the course of this discussion, he proved the following result. 

\REFTHM{main}(W. Thurston, December 2010)\label{billthm} 
 An oriented planar graph $\G$  with  $2d-2$ vertices of valence $4$
 is equal to $f^{-1}(\Sigma)$ for some degree $d$ branched cover $f:S^2\to S^2$ and   some $\Sigma$ if and only if 
\begin{enumerate}
\item every face of $\G$ is a Jordan domain (in particular $\G$ is connected),
\item  (global balance) for any alternating white-blue coloring of the faces of $\G$, there are $d$ white faces, and there are $d$ blue faces, and 
\item (local balance) for any oriented simple closed curve in $\G$ that keeps the blue faces on the left and white on the right (except at the corners), there are strictly
more blue faces than white faces on the left side.
\end{enumerate}
\ENDTHM

\begin{figure}[h]
   \centering
   \includegraphics[width=4in]{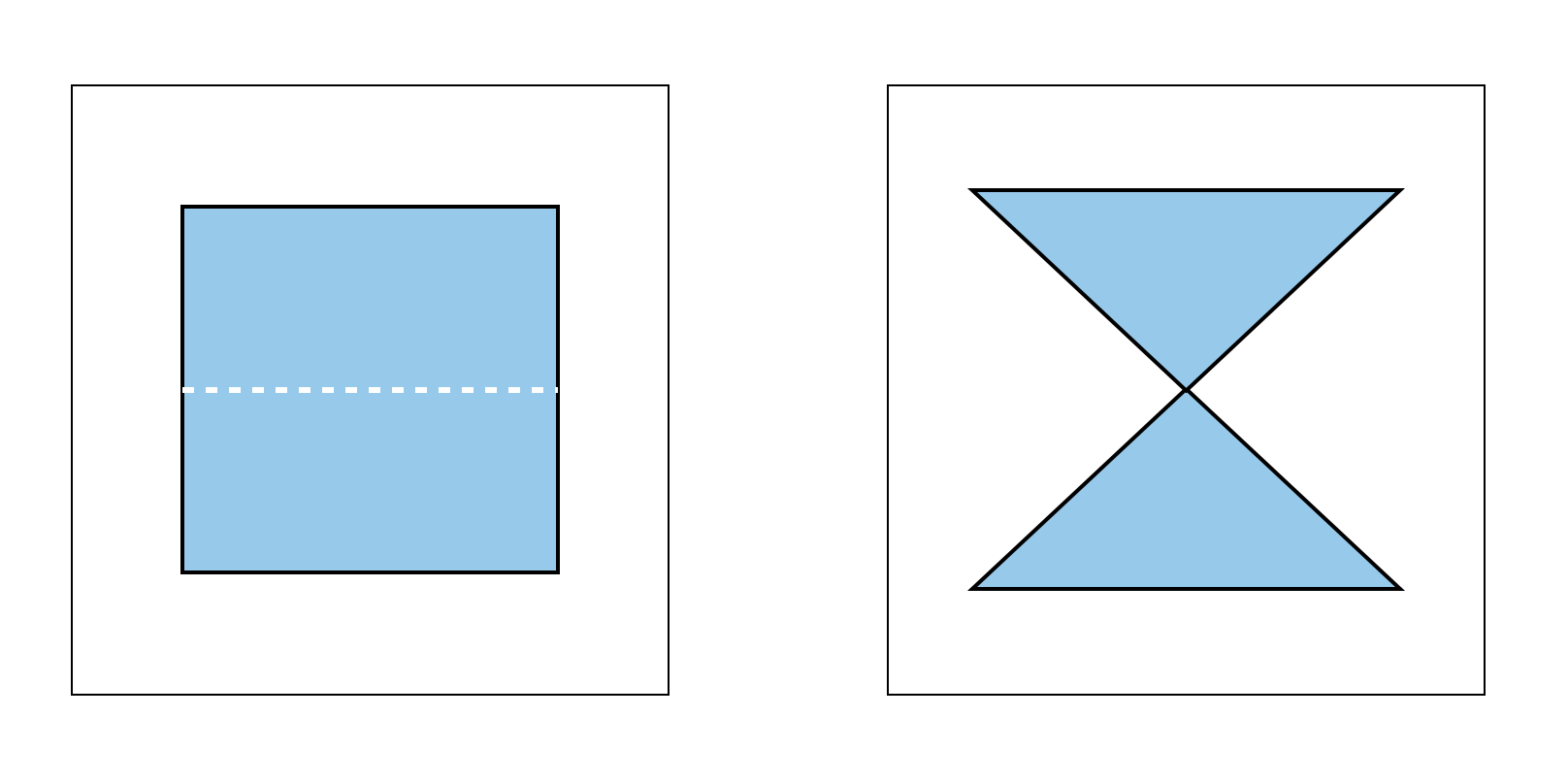}
   \caption{\small If the dashed curve in the blue region is pinched to a point, the number of blue faces increases by $1$, while the number of white faces remains constant.}
   \label{pinch} 
\end{figure}
\begin{figure}[htbp] 
   \centering
   \includegraphics[width=3.5in]{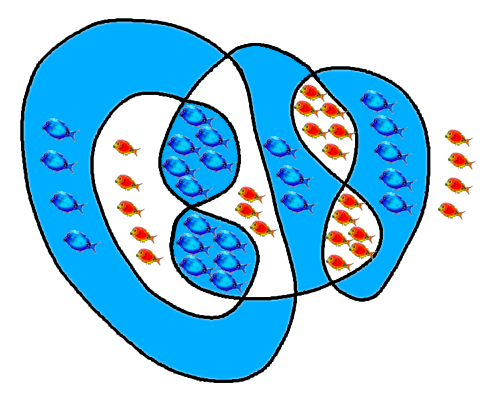} 
   \caption{Give the surface an orientation, and orient the graph above. There are equal numbers of blue and white faces, so the graph $\Gamma$ is globally balanced; however, the right side of the graph violates the local balance condition, regardless of the orientation chosen. Note that in each face of this diagram, the number of  fish plus the number of corners is equal to $8=2(5)-2$. 
This figure is a slightly modified version of one made by W. Thurston.}
   \label{Localimbalanced}
\end{figure}

 Any oriented $4$-valent graph  satisfying the three conditions above is called {\em balanced}. One way to construct arbitrary 4-valent graphs which  fail the global balance condition in Theorem \ref{billthm} is to pinch arcs in regions of one color only. The number of faces of that color increases without changing the number of the other faces, see Figure \ref{pinch}. 
A globally balanced but locally imbalanced graph is shown in Figure \ref{Localimbalanced}. Notice that the balance conditions are insensitive to the orientation of the surface, and they are insensitive to the choice of blue and white faces in an alternating coloring of the faces.

\begin{figure}[h] 
   \centering
   \includegraphics[width=2.5in]{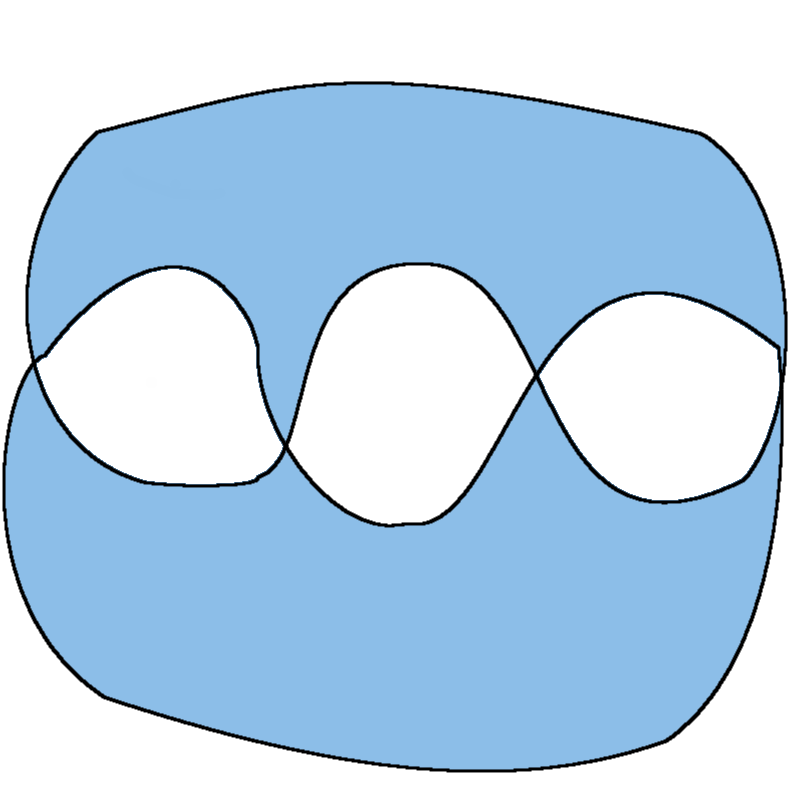} 
   \caption{A globally imbalanced graph. There are $4$ white faces against $2$ blue ones. This figure was made by W. Thurston.}
   \label{Imbalanced}
\end{figure}

%

  Before giving the proof of Theorem \ref{billthm}, we will play with some examples and study the graphs arising from maps $f:S^2\to S^2$ of degree $4$ (which is exactly how Bill invited us to participate in the original email discussion where this work began). 
 
  \noindent {\bf Outline.} In Section \ref{examples}, we begin by discussing Hurwitz numbers for generic branched covers $S^2\to S^2$ of degree $d$, and we then present a complete list of graphs arising in the case $d=4$. We then establish the proof of Theorem \ref{billthm} in Section \ref{proof}, closely following Thurston's original proof. We have added some details (mainly for our own understanding). In Section \ref{decomposition}, we present some of W. Thurston's work on decomposing balanced planar graphs into `standard pieces'. In the appendix, we prove Hurwitz's classical result, that (up to isomorphism) the number of generic degree $d$ branched covers $S^2\to S^2$ is 
\[
\frac{(2d-2)!d^{d-3}}{d!},
\]
presenting a geometric proof due to Duchi-Poulalhon-Schaeffer in \cite{DPS}. 
 
  \noindent {\bf Acknowledgments.} First and foremost, we would like to thank Bill Thurston for starting this email discussion in 2010, engaging a large number of mathematicians in this project, and sharing his work with all of us. It is with great pleasure that we present this part of his work. Thanks to Laurent Bartholdi, Dylan Thurston, and Jerome Tomasini for helping with this manuscript in various ways. The research of Sarah Koch is supported in part by NSF grant DMS 1300315 and a Sloan Research Fellowship, and the research of Tan Lei is supported in part by ANR LAMBDA. 
 
\section{Examples and a complete list of quartic diagrams}\label{examples}

\subsection {Hurwitz numbers.} Let $f:S^2\to S^2$ be a generic orientation-preserving  branched cover of degree $d\geq 2$.
Label the critical values $R_f=\{v_1,\ldots,v_{2d-2}\}$, and let $\Sigma$ be an oriented Jordan curve running through the critical values following the order of the labels. 
Color the left complementary disc of $\Sigma$ blue, and color the right complementary disk white. 
Consider $\Sigma$ as a graph with $V_\Sigma=R_f$, where each $v\in V_\Sigma$ has valence $2$. 

  Then $f^{-1}(\Sigma)$ is an oriented graph with vertices $V_{f^{-1}(\Sigma)}=f^{-1}(R_f)$; there are precisely $2d-2$ vertices in $V_{f^{-1}(\Sigma)}$ with valence $4$ (corresponding to the critical points of $f$), and the other vertices in $V_{f^{-1}(\Sigma)}$ have valence $2$. We label the vertices of $f^{-1}(\Sigma)$ as follows: for each $c\in V_{f^{-1}(\Sigma)}$,  label the vertex $c$ with the number $j$ if $f(c)=v_j$. In this way, there are $d-1$ vertices in $V_{f^{-1}(\Sigma)}$ labeled $1$, $d-1$ vertices in $V_{f^{-1}(\Sigma)}$ labeled $2$, and so on. Note that because $f:S^2\to S^2$ is generic, for each $1\leq j\leq 2d-2$, there is a unique vertex labeled $j$ with valence $4$. Color each face  of $f^{-1}(\Sigma)$ either white or blue according to the color of its image, under $f$, in $S^2-\Sigma$. The result is a checkerboard pattern in the domain of $f:S^2\to S^2$. 

  Two such labeled graphs  are said to be {\em isomorphic} if there is an orientation-preserving homeomorphism
$\phi:S^2\to S^2$ sending one oriented graph to the other, preserving the labels of the vertices. According to Hurwitz, the number of isomorphic labeled  graphs (of degree $d$) is 
\[
\frac{(2d-2)!}{d!} d^{d-3}
\]
(see the appendix for further details). When $d=4$, this number is $120$.

\subsection{The case $d=4$: cataloguing all 120 quartic diagrams} We illustrate the complete list of all 120 labeled graphs arising from  $f:S^2\to S^2$, where $d=4$. They are organized into groups, catalogued by their underlying $4$-valent graphs.  The data are are based on a 
list, computed by L. Bartholdi with Gap, of $6$-tuples of transpositions in $S_4$, with trivial
product, modulo diagonal conjugation by $S_4$. The number of such tuples is equal to the
number of isomorphism classes of generic degree $4$ maps $f:S^2\to S^2$. 

\subsubsection{Constructing $\Gamma$ from $f:S^2\to S^2$} In the case $d=4$, we recall the steps used to construct $\G$, the oriented underlying 4-valent graph of the branched cover $f:S^2\to S^2$: 
\begin{enumerate}
\item Take an oriented Jordan curve $\Sigma\subseteq S^2$ in the range running through the 6 critical values $R_f=\{v_1,\ldots,v_6\}$, and consider $f^{-1}(\Sigma)\subseteq S^2$ as an oriented graph in the domain. 
\item There are $6$ vertices in $f^{-1}(\Sigma)$ with valence 4 (corresponding to the critical points of $f$), there are 12 vertices in $f^{-1}(\Sigma)$ with valence 2 (the cocritical points). 
\item Each face of $f^{-1}(\Sigma)$ has $6$ vertices on its boundary.
\item The underlying oriented 4-valent graph $\G$ is the oriented graph $f^{-1}(\Sigma)$ obtained by forgetting the vertices which have valence $2$. 
\end{enumerate}
\subsubsection{Constructing $f:S^2\to S^2$ from a diagram}
\begin{figure}[h]
   \centering
   \includegraphics[width=5in]{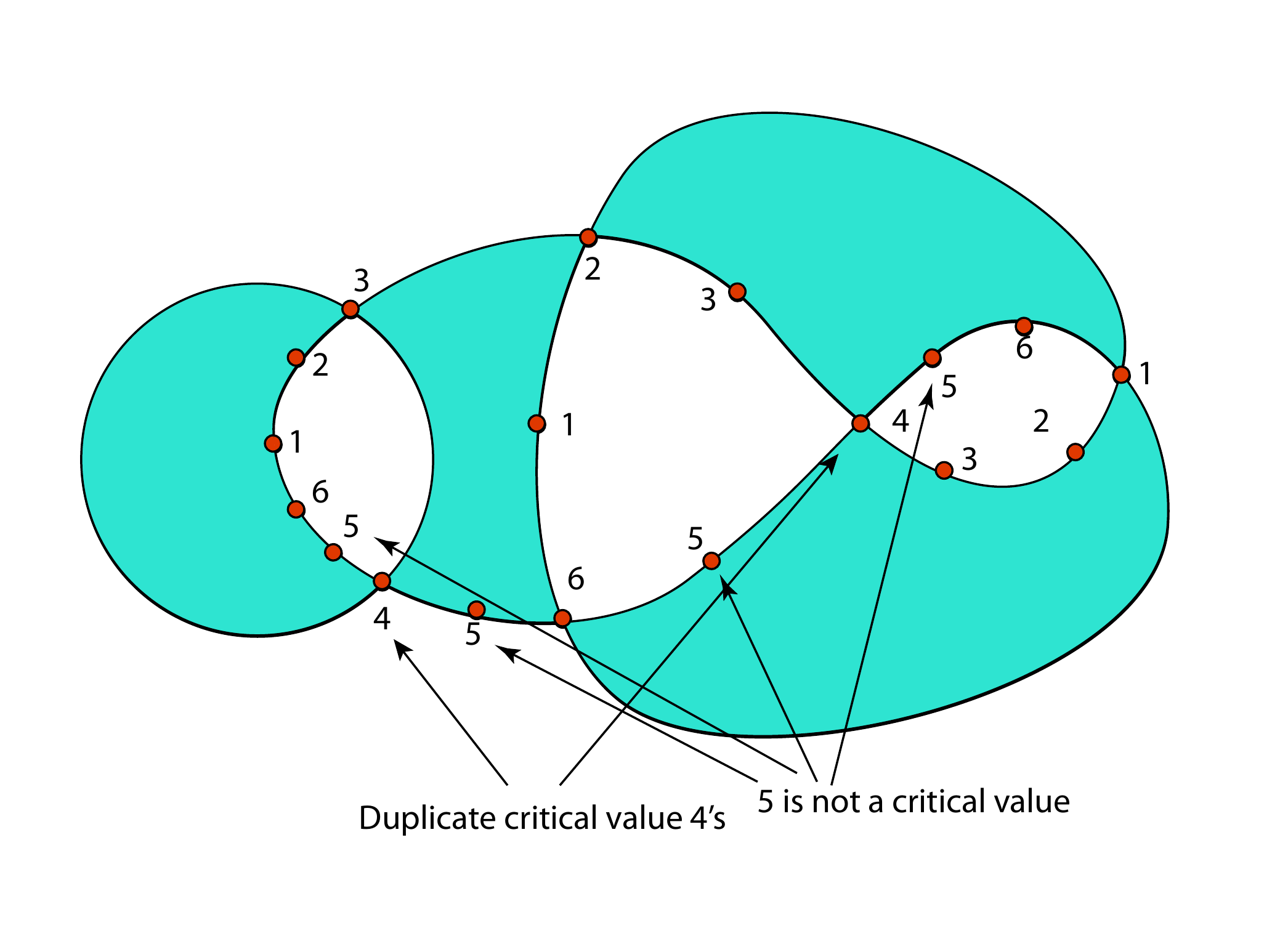} 
   \caption{Note how the numbering goes clockwise around each white face and counterclockwise around each blue face. There's a duplicate critical value:  the two vertices labeled 4 map to the same point in the codomain, and none of the vertices labeled 5 is a critical point, so there are only five distinct critical values in the codomain, corresponding to the labels: $1,2,3,4,6$. This figure was made by W. Thurston.}
   \label{duplicate}
\end{figure}

For each of the diagrams below, we build a degree $4$ branched cover $f:S^2\to S^2$ in the following way:
\begin{enumerate}
\item  Choose an alternating coloring of the faces (a checkerboard-coloring). 
\item Choose one 4-valent vertex (that is, a critical point) to map to the critical value $v_1$; label this critical point with a `$1$'. 
\item Label the remaining cocritical points and critical points so that they appear in counterclockwise order on the
boundary of each blue face, see Figure \ref{duplicate} for this labeling scheme. 
\item Define the covering $f$ first on the vertices  following the labeling, and then on the edges by homeomorphic extension, and finally on the faces by homeomorphic extension (with a little help from the Sch\"oenflies theorem). 
\end{enumerate}

\begin{figure}[h] 
   \centering
   \includegraphics[width=3in]{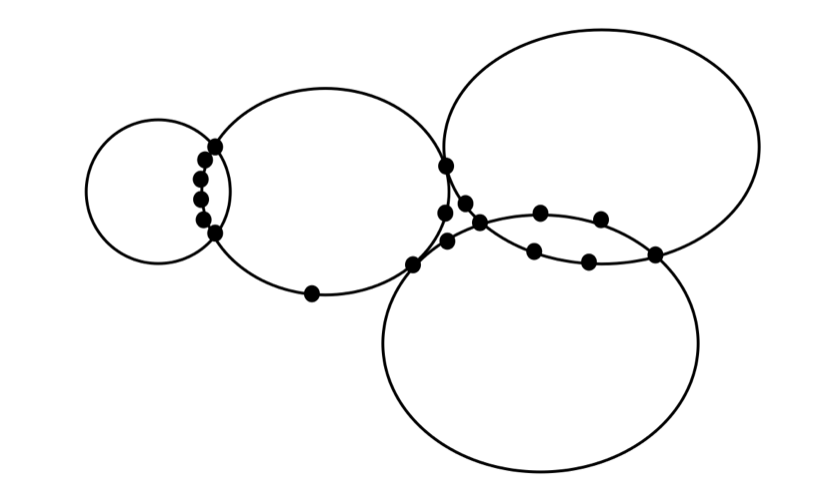} 
   \caption{This diagram cannot be realized as the underlying diagram for any $f:S^2\to S^2$. In fact, any labeling of this diagram will result in a `duplicate critical value.' See Figure \ref{duplicate}.}
   \label{badpic}
\end{figure}

\begin{remark}\label{notreal}
The recipe above for building a map $f:S^2\to S^2$ from the labeled diagram does not always work; for example, there is a labeled diagram that satisfies the conditions above, but there is no branched cover $f:S^2\to S^2$ which gives rise to the labeled graph in Figure \ref{duplicate}, or in fact to the underlying unlabeled graph in Figure \ref{badpic}. 
\end{remark}

  Once a diagram is realizable as an underlying oriented $4$-valent graph for some $f$ and $\Sigma$, there are 
\begin{itemize}
\item 6 choices of the critical point to label with `$1$', and 
\item 2 choices of the set of white/blue faces.
\end{itemize}
So in general, each diagram below gives rise to 12 distinct labeled diagrams. However, some of these diagrams have extra symmetry, and in these cases, there are fewer than 12 labeled diagrams (up to isomorphism). 


\begin{figure}[h] 
   \centering
   \includegraphics[width=6in]{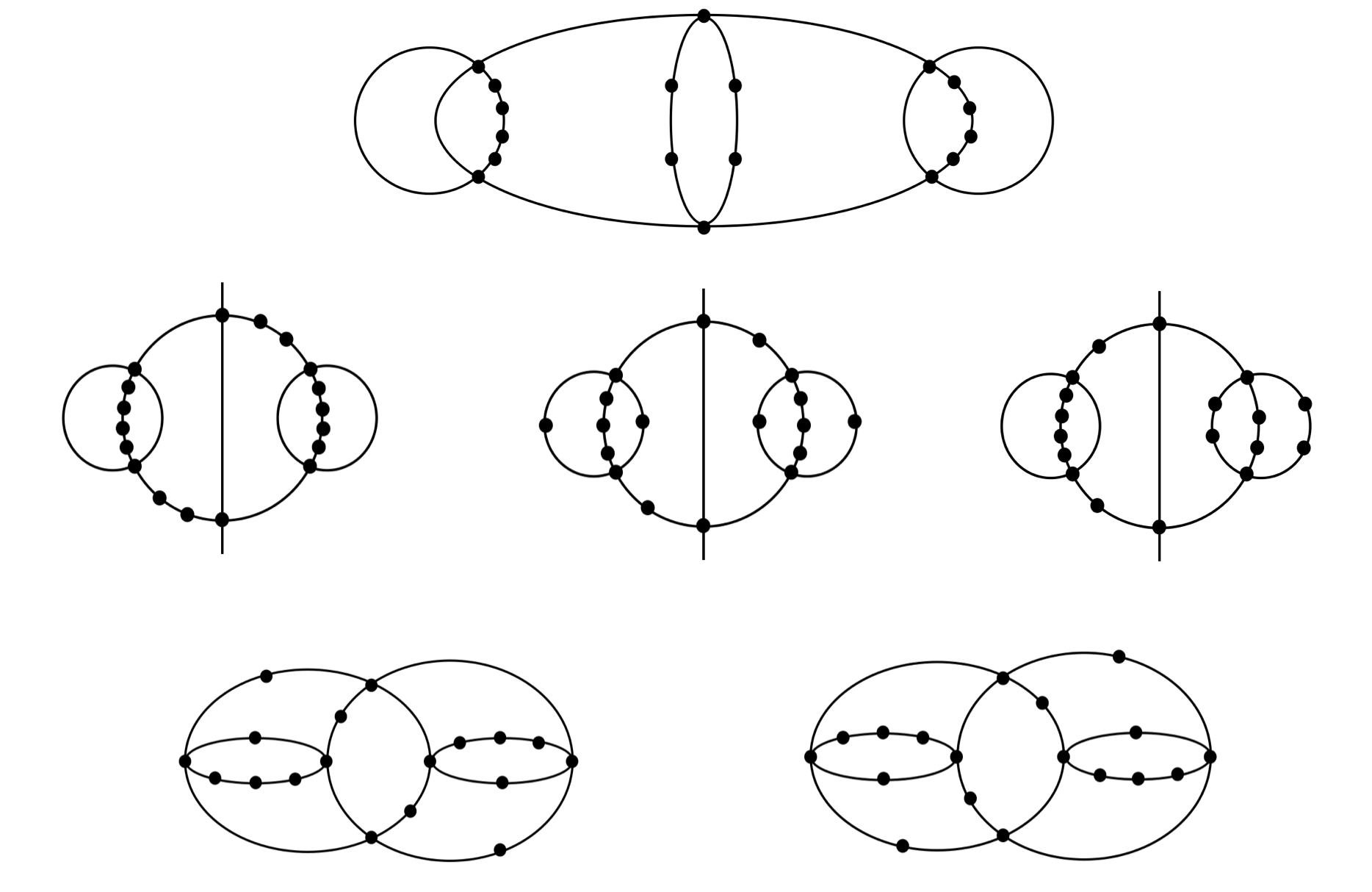} 
   \caption{These six graphs represent 36 branched covers; each graph above gives rise to three distinct branched covers $f:S^2\to S^2$, up to isomorphism.}
   \label{quarticgraph2}
\end{figure}

\begin{figure}[h] 
   \centering
   \includegraphics[width=6in]{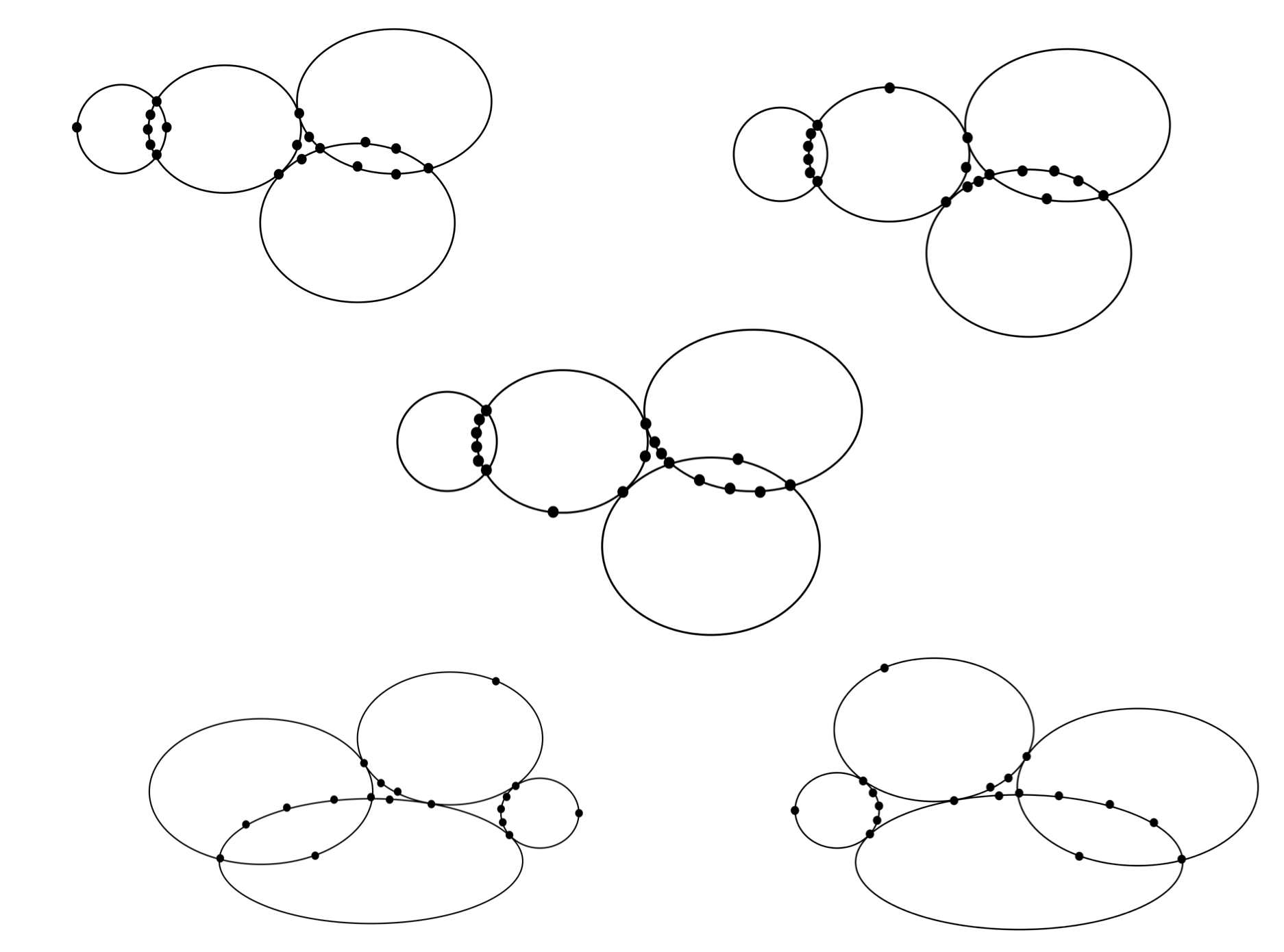} 
   \caption{These five graphs represent 60 branched covers; each graph above gives rise to $12$ distinct branched covers $f:S^2\to S^2$, up to isomorphism.}
   \label{fig:example}
\end{figure}

\begin{figure}[h] 
   \centering
   \includegraphics[width=6in]{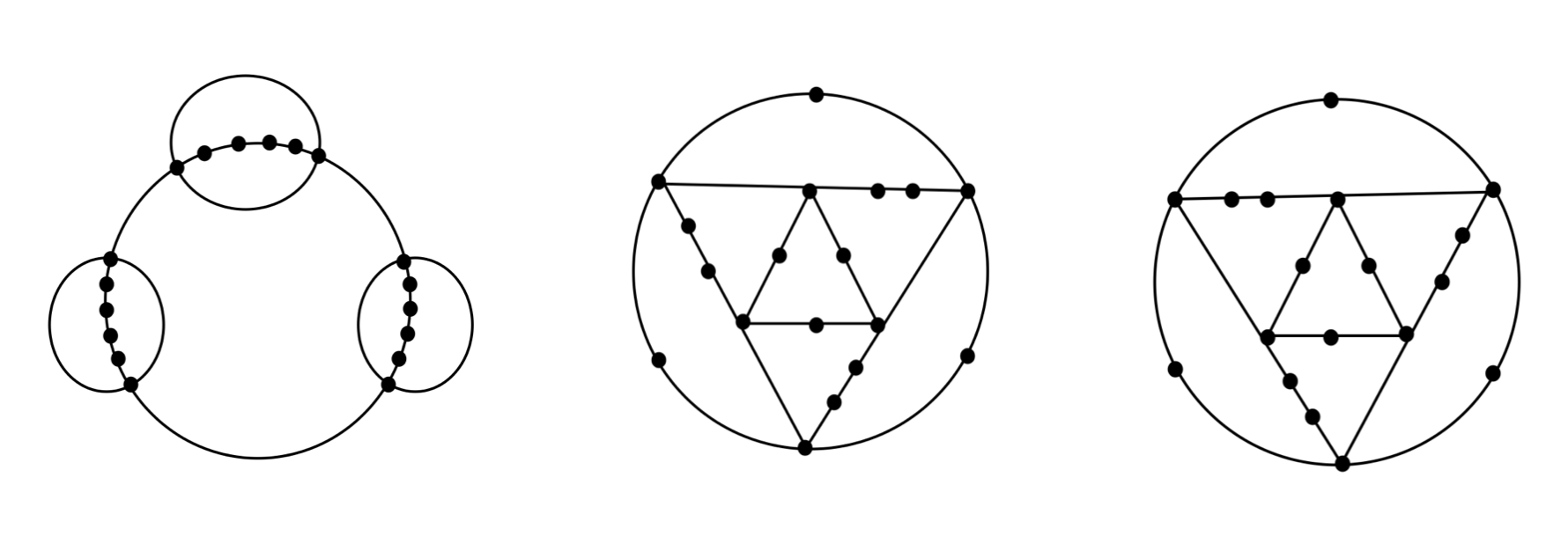} 
   \caption{These three graphs represent 6 branched covers; each graph above gives rise to two distinct branched covers $f:S^2\to S^2$, up to isomorphism.}
   \label{bye}
\end{figure}

\begin{figure}[h] 
   \centering
   \includegraphics[width=4in]{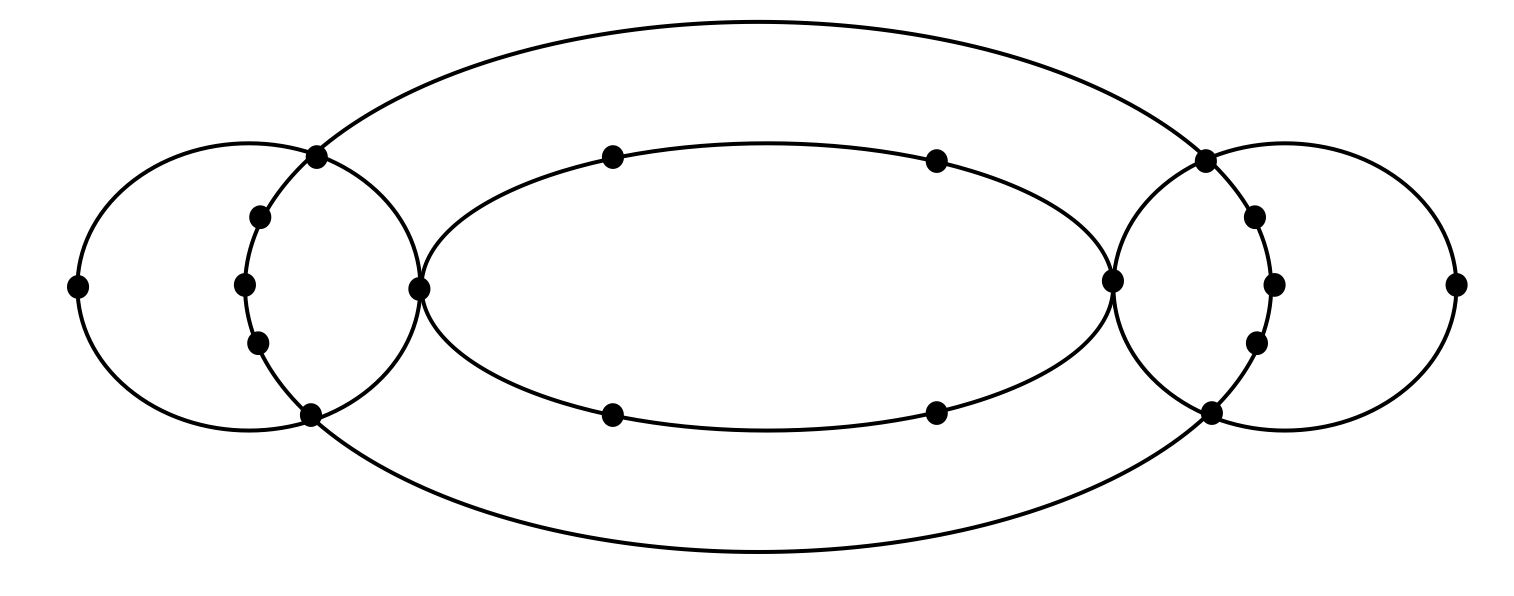} 
   \caption{This graph represents 6 branched covers up to isomorphism.}
   \label{hi}
\end{figure}

\begin{figure}[h] 
   \centering
   \includegraphics[width=6in]{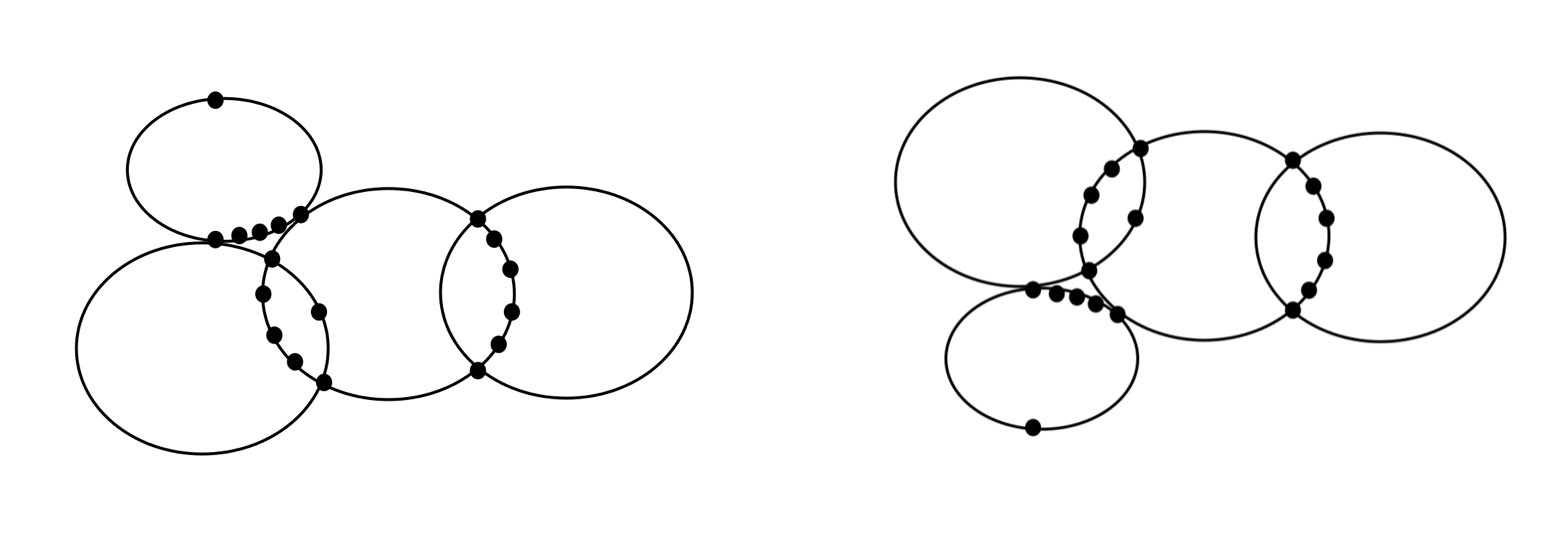} 
   \caption{These two graphs represent 12 branched covers; each graph above gives rise to six distinct branched covers $f:S^2\to S^2$, up to isomorphism.}
   \label{quarticgraph}
\end{figure}


\section{Proof}\label{proof}

  We will prove Theorem \ref{billthm}, following the main ideas in W. Thurston's original proof. Recall the statement:

  {\em An oriented planar graph $\G$  with  $2d-2$ vertices  of valence $4$
 is equal to $f^{-1}(\Sigma)$ for some degree $d$ branched cover $f:S^2\to S^2$ and some $\Sigma$, if and only if 
\begin{enumerate}
\item every face of $\G$ is a Jordan domain (in particular $\G$ is connected), 
\item  (global balance) for any alternating white-blue coloring of the  faces of $\Gamma$, there are $d$ white faces, and there are $d$ blue faces, and 
\item (local balance) for any oriented simple closed curve in $\G$ that keeps the blue faces on the left and white on the right (except at the corners), there are strictly
more blue faces than white faces on the left side.
\end{enumerate}
}
 \beginp The key idea is to first translate the realization problem into finding a pattern of vertices (either 2-valent or 4-valent) so that each face of $\Gamma$ has precisely $n:=2d-2$ vertices on the boundary, and then to reduce the problem into a matching problem in graph theory. Indeed, let $\Sigma$ be an oriented simple closed curve which goes through the critical values of a degree $d$ branched cover $f:S^2\to S^2$. Consider the oriented graph $f^{-1}(\Sigma)$; it has $n$ vertices of valence $4$, and $n(d-1)-n$ vertices of valence $2$. 
Split each 2-valent vertex into two `dots', putting one into each of the two adjacent neighboring  regions. Then split each 4-valent vertex (or corner) into $4$ `dots', putting one into each of the four neighboring regions. Every face of $f^{-1}(\Sigma)$ then contains precisely $n$ dots in its interior. Grouping the dots back together at the vertices then becomes a matching problem. 

  Now take an arbitrary oriented graph $\G$ with $n$ 4-valent vertices and no other vertices. Color the faces in a blue and white checkerboard pattern.
Put $n$ men in each blue  face (represented by blue fish in Figure \ref{Localimbalanced}), and $n$ women in each white face (represented by red fish in Figure \ref{Localimbalanced}). Each is trying to find a partner from one of the neighboring faces.   In addition, each corner requires one person per neighboring face (that makes two women and two men).     Is there a perfect matching?   The usual marriage criterion in the marriage theorem from graph theory states that for every set of $N$ women, there are at least $N$ men who are potential mates; this can be reduced to the global balance and local balance conditions in the theorem. Here are the details:

\medskip

\noindent {\bf Necessity, condition 1:} complementary Jordan domains.

 \noindent Let $f:S^2\to S^2$ be a generic orientation-preserving branched cover of degree $d$, and let $\Sigma$ be an oriented simple closed curve passing though the $n$ critical values of $f$. Note that $f$ restricts to a degree $d$ covering map $S^2-f^{-1}(\Sigma) \to S^2-\Sigma$; it follows that the connected components of $S^2-f^{-1}(\Sigma)$ are Jordan domains. 

\medskip

\noindent {\bf Necessity, condition 2:} global balance.

 \noindent Change coordinates in the range so that  $\Sigma$ becomes the unit circle $S^1$, and orient it so the unit disk is on the left, and color the unit disk blue. The white face is the complement. We may also assume that the point $1$ is not a critical value. Now the oriented graph we are considering in the domain is $f^{-1}(S^1)$. 
Each of the blue faces in $S^2-f^{-1}(S^1)$ contains one zero of $f$, and each of the white faces in $S^2-f^{-1}(S^1)$  contains one pole of $f$.  The global balance condition is equivalent to the condition that the number of zeros of $f$ equals the number of poles of $f$. 

\medskip

\noindent {\bf Necessity, condition 3:} local balance.

 \noindent Choose an alternating coloring of the complementary components of $f^{-1}(\Sigma)$. 
Let $\La$ be a Jordan domain bounded by an oriented curve $\g$ in $f^{-1}(\Sigma)$ keeping blue faces on the left.
Suppose that in $\La$ there are $W$ white faces, $B$ blue faces, $E$ edges in the interior of $\La$, and $X$  edges on the boundary.
The total number of women in the white faces is  $nW - E$.
The total number of men in the blue faces is $nB - E - X$, with $X>0$.
For these $nW - E$ women, all of their potential mates are among the $nB - E - X$ men. Necessary and sufficient conditions from the marriage lemma imply $nW - E\le nB - E - X$, from whence it follows that $B> W$. This establishes the local balance condition of the theorem.


\medskip

\noindent {\bf Sufficiency, step I:} adding $2$-valent vertices to $\G$.

  \noindent Suppose that $\Gamma$ is an oriented $4$-valent graph which satisfies the conditions of the theorem. Every face $F$ of $\G$ is a Jordan domain, so the number of edges on the boundary of $F$ is equal to the number of corners of $F$. 
For any face with $k$ corners, we have  $k\le n$. After matching the corners, there are $n-k$ people left in the face to match. Since  the numbers of faces of the two colors are equal, the total number of men who are not matched with corners is equal to the total number of women who are not matched with corners. 
 
  In what follows, we will  consider the men and women who have not been matched with corners.  Let $S$ be any subset of these women. We want to show that the set of potential mates of $S$ is at least as large as $S$ (satisfying the condition in the marriage theorem).
 
  If we augment the set $S$ with all women in every face containing an element of $S$, this will not increase the number of their potential mates, so we may assume  $S$ is actually the union of all women in some collection $U$ of white faces.  Let $R$ be the closure of  $U$ with all adjacent neighboring blue faces of $U$. Then the men in $R$ are exactly the potential mates for the women in $S$.  The orientation of the boundary edges of $R$ leaves $R$ on the left, and leaves blue faces on the left as well.

  If the interior of $R$ is not connected, then the women in one connected component can only match with men in that same connected component.
We  will therefore establish that there are enough mates  for the elements of $S$ in  each connected component of the interior of $R$, from whence it will follow that there are enough  mates in $R$ for all of the elements of $S$. To this end, we assume that the interior of $R$ is connected, so that every complementary component of $R$ in $S^2$ is a Jordan domain.

  Denote the number of men/women inside/outside $R$ by the following tableau:

\begin{center}\begin{tabular}{c|c|c} & number of men \mystrut & number of women \\ \hline
inside $R$ & $m$  \mystrut & $w:=|S|$ \\ \hline
outside $R$ & $m^*$ \mystrut & $w^*$
\end{tabular}\end{center}

  Assume at first that $R$ is simply-connected, that is, $R$ is bounded by a simple closed curve leaving blue faces on the left.

  Denote by:
\begin{itemize}
\item $X_1$ the number of boundary  corners of $R$ which are on the boundary of only one face in $R$ (this face is necessarily blue);\item
 $X_3$ the number of boundary  corners of $R$ that are on the boundary of three faces in $R$ (necessarily  two blues and one white);\item
 $v$ the number of vertices  in the interior of $R$;\item
$W$ the number of white faces in $R$ and\item $B$ the number of blue faces in $R$.
\end{itemize}
Note that every boundary corner  of $R$ belongs to one of the  two types either  $X_1$ or  $X_3$, and
$$w=|S|=	nW -(2v+ X_3),\quad m=n B -(X_1+ 2X_3+2v).$$

  By the local balance condition, we have
$W \le B-1$.  Therefore,  using  $X_1+ X_3 \le n$, we get
\begin{eqnarray*} w=nW -(2v+ X_3) & \le&  n ( B-1)- (X_3+ 2v) \\ &=&  n B -(n+ X_3+2v)\\   &\le &  n B -(X_1+ 2X_3+2v)=m.
\end{eqnarray*}
This is the desired inequality.

  We now deal with the case where $R$ is not simply-connected. The number of women in the complement of $R$ is $w^*$, and the number of men in the complement of $R$ is $m^*$. If we can show that $m^*\leq w^*$, the fact that $m+ m^*=w + w^*$ will imply that $w\le m$, which is the desired inequality. 

  Let $Y$ be a complementary component of $R$. It is bounded by a Jordan curve keeping $Y$ and white faces on the left.
The same curve in the opposite direction, enclosing the complement of $Y$,
  has blue faces on the left. So by the local balance condition there are more blue faces than white faces outside $Y$. Using the assumption that there are equal numbers of blue and white faces (the global balance condition), there are more white faces than blue faces in $Y$.
We then apply the same argument above in the case that $R$ is simply-connected (replacing $R$ by $Y$ and reversing the colors) to conclude
that the set of men in $Y$ have enough mates (necessarily in $Y$); this establishes that $m^*\leq w^*$. 

  We have now proved that for any set $S$ of women, the number of potential mates is at least as large as $S$. By the marriage theorem, there is a perfect matching. For each matched pair, merge them  into a vertex on one common boundary edge of the two faces. Make sure that
the  vertices are pairwise disjoint. We have now enriched $\G$ into a graph with $n(d-1)$ vertices, where $n$ are 4-valent and the rest are 2-valent.
Furthermore, every face has exactly $n$ vertices on its boundary.

\medskip

\noindent{ \bf Sufficiency, step II:} consistent labeling of the vertices.  

  \noindent We now want to label the enriched vertices of $\G$ by the set $\{1, \ldots, n\}$, such that each label appears exactly $d-1$ times, and the natural order of the labels appears counterclockwise around each blue face (and therefore clockwise around each white face), as shown in Figure \ref{duplicate}.

  Consider $V_\G\times \Z/n\Z$, where $V_\G$ denotes the set of vertices of $\Gamma$.
For every edge $e$,  from vertex $v$ to vertex $v'$ 
define $h_e : \{v\}\times \Z/n\Z\to \{v'\}\times \Z/n\Z$ by 
\[
h_e:(v,x)\mapsto (v',x+1\;\mathrm{mod}\; n).
\]
The collection of $h_e$ forms a
cocycle with coefficients in $\Z/n\Z$. Since $H^1(S^2)=0$, this is the coboundary of a function to $\Z/n\Z$.  
 
  More generally: consider a cocycle with values in $S^1 = \R/2 \pi \Z$, and parametrize $\R \mathbb{P}^1$ by a circle-valued coordinate, then a particular cocycle yields exact positions for the critical values, up to rotations of the circle (constant of integration).

\medskip
 
\noindent { \bf Sufficiency, step III:} avoiding duplicate critical values.

  \noindent It may happen, has we see in the quartic diagrams (see Figure \ref{duplicate}) that  two critical points have the same label. In this case, we may perturb the map to obtain matchings and labels so that the $n$ critical points have pairwise distinct labeling.
 \qed

 \subsection{Further comments} 
 
\noindent It would be nice to have a better proof of  sufficiency that can be turned  into an algorithm to construct arbitrary rational functions whose critical values are on the unit circle. There should be some kind of cell structure whose lower-dimensional cells correspond to subdivisions of the sphere with vertices of higher even order as critical points coalesce, obtained by collapsing forests contained in the 1-skeleton of the subdivisions in the generic case. Collisions between critical points is related to the strata in Hurwitz's study.

\medskip

  \noindent {\bf Graph flow.} One way to test whether a graph is balanced or not is to directly test a possible matching.  There are good algorithms for finding matchings when they exist, or finding an obstruction when they don't exist: this is very related to graph flow, which has fast algorithms. 

  Here is a graph flow version of the balance condition: let $\G$ be a connected 4-valent graph with an alternating 
coloring of its faces. Give each face $F$ a weight $w(F)$ equal to the number of vertices of $\G$ minus the number of corners of $F$. Define an abstract graph flow $\Sigma$  whose vertices are $\{v(F), \ F \text{ a face of }\G\}$ and whose directed edges are of the form $\{(v(F), v(F')\}$  for every ordered pair $(F, F')$ of blue-white faces sharing one or more boundary edges in $\G$. 
Add two vertices $D$ (for departure) and $A$ (for arrival) to $\Sigma$. Direct an edge from $D$ to every  vertex $v(F)$ of blue face $F$. This edge has  capacity $w(F)$. Direct an edge from every   vertex $v(F')$ of white face $F'$ to $A$, with capacity $w(F')$.  Now the map is balanced if the maximal flow reaches the maximal capacity $\sum_{F\text{ blue}} w(F)$.

\section{Decompositions of balanced planar graphs}\label{decomposition}

  Our objective in this section is to understand the structure of balanced graphs from the point of view of decomposing them into standard pieces. 

\subsection{A *22 decomposition}

\begin{figure}[h] 
   \centering
   \includegraphics[width=5in]{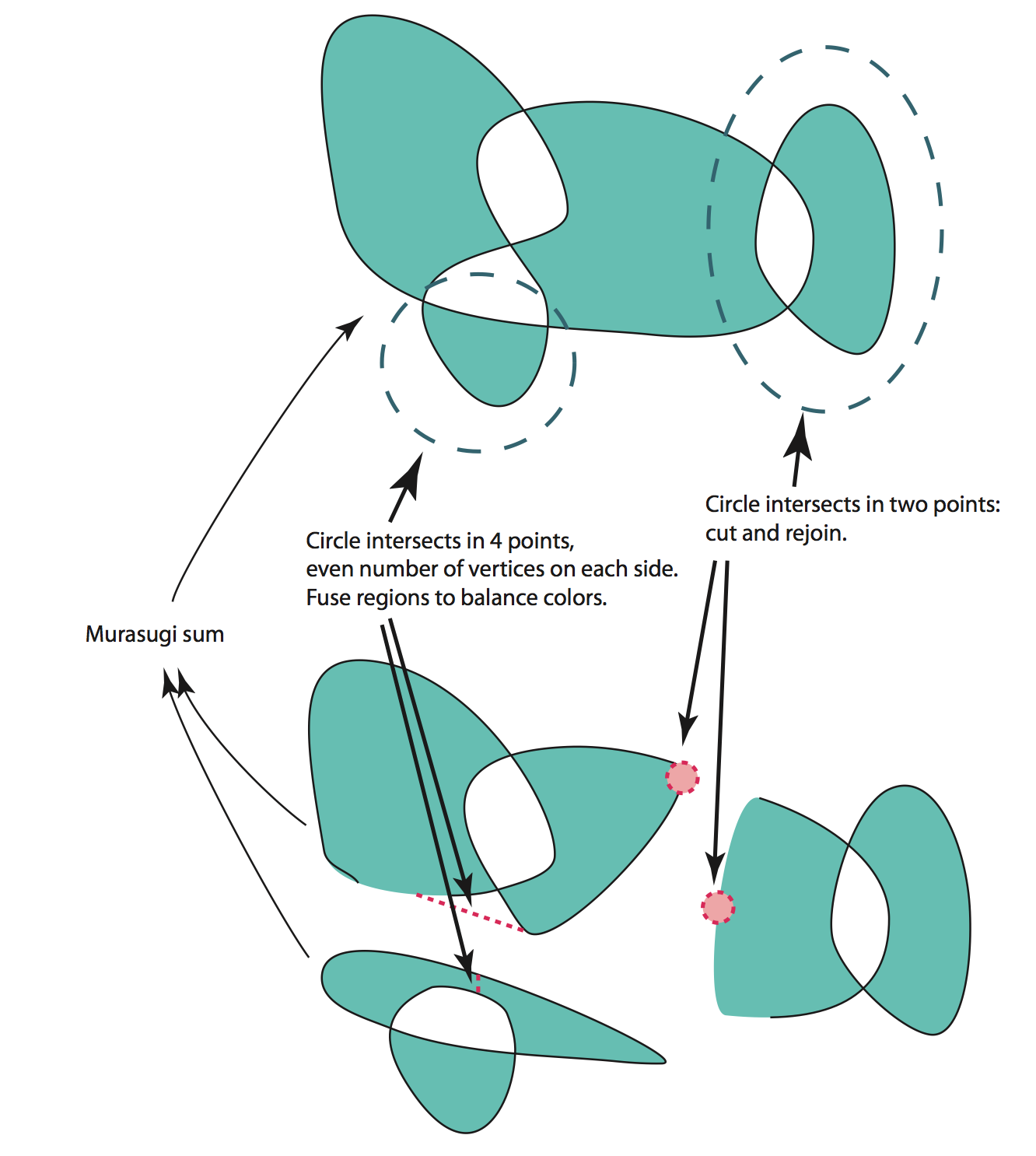} 
   \caption{The dashed Jordan curve on the right cuts $\G$  in exactly two points and they are on distinct edges; it is nontrivial. This figure was made by W. Thurston.}
   \label{22-decomposition}
\end{figure}

%

  Start with a planar $4$-valent graph $\G$. Find a nontrivial Jordan curve $\g$ in $S^2$ that intersects $\G$ in only two points
(here nontrivial means that the two intersection points
are on different edges). This happens exactly when 
 two faces of opposite colors have more than two common edges. Cut  $S^2$ along $\g$ and collapse each of the wound curves into a single point (these points will become regular points of the surgered graph).
Continue this process on the two new diagrams until there are no longer any nontrivial curves with two intersection points.

\subsection{Two types of tangle decompositions.}

\begin{figure}[htbp] 
   \centering
   \includegraphics[width=5in]{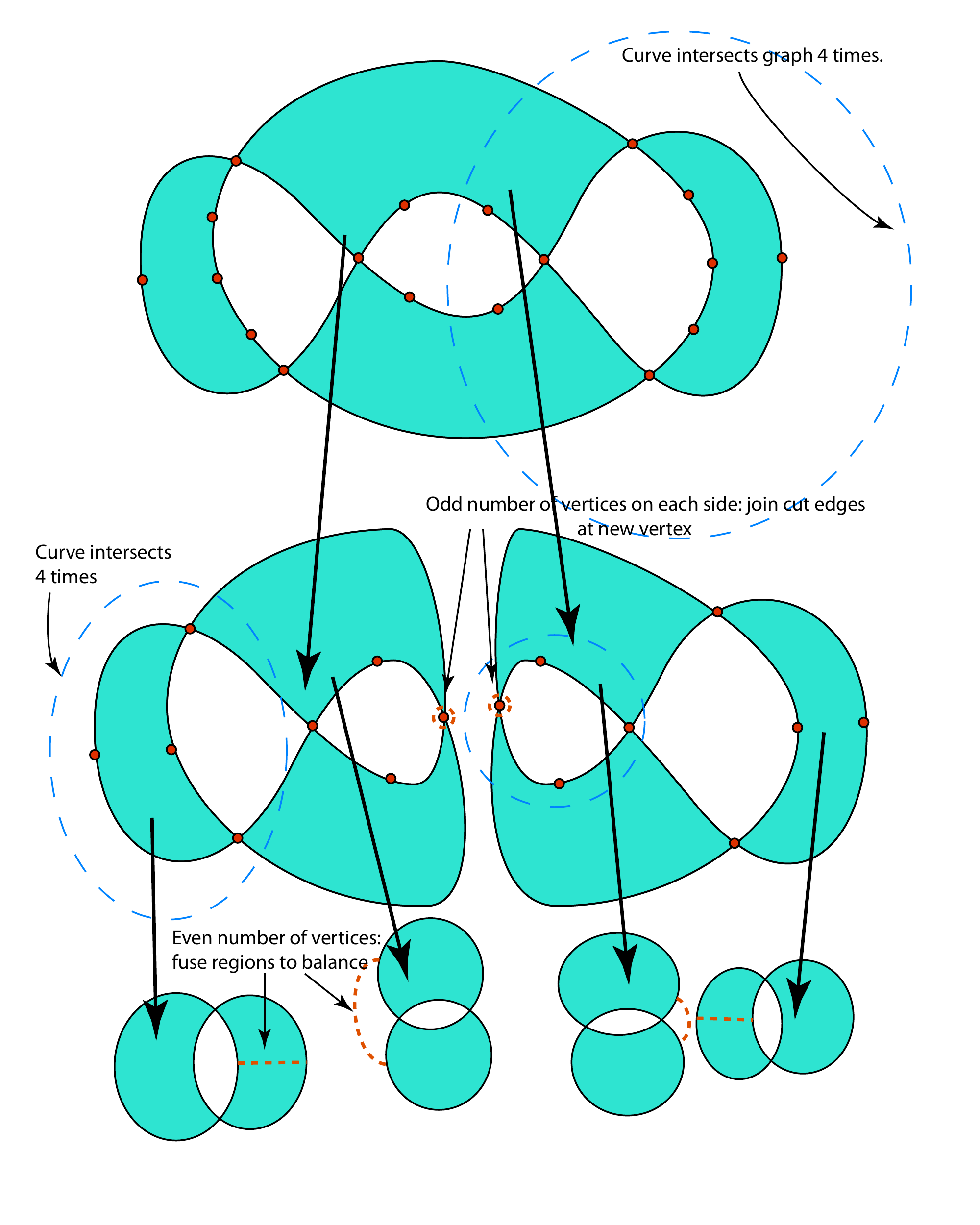} 
   \caption{The two types of tangle decompositions, the 2-valent vertices are not taken into consideration, they should be added later. This figure was made by W. Thurston.}
   \label{tangle}
\end{figure}


  Now look for a nontrivial Jordan curve  that intersects $\G$ in four points, here nontrivial means that the curve does not go around a single vertex,  and that the four points are on four distinct edges.

  If the curve $\g$ separates the $4$-valent vertices of $\G$ into a pair of odd numbers of vertices, 
cut  $S^2$ along $\g$ and collapse each of the wound curves into a single point. The single point becomes a new $4$-valent vertex of the surgered graph.
See the top right dashed curve in Figure \ref{tangle}.

  Assume that the curve $\g$ separates the vertices of $\G$ into a pair of even numbers of vertices, and assume that $\G$ has an equal number
of white and blue faces. Cut the sphere $S^2$ along $\g$. Then the color imbalance on the two sides of the cut is $\pm 1$. 
For the side with more white (half)-faces, glue the two white edges of the cut copy of  $\g$ to connect the two white half faces into one singe white face
(and fold the two blue edges into two segments),
and symmetrically for the other side with more blue (half)-faces. See for example the middle two dashed curves in Figure \ref{tangle} and 
the top left dashed curve in Figure \ref{22-decomposition}.

  The converse procedure of the last decomposition corresponds to {\em Murasugi sum} of knots. 
To do Murasugi sum of two diagrams, consider them as embedded in two spheres. Place one sphere on the left and one on the right (in $3$-space). Remove a rectangle from oppositely-colored regions of the two  diagrams,  
where the rectangles have two edges on different sides of  a region and two edges interior to the region.  
Then glue them so as to match colors.  The two spheres are merged into a single one.

\subsection{Indecomposable pieces}

Every diagram decomposes by these three operations into 
\begin{itemize}
\item pieces that are isomorphic the unique quadratic diagram ($2$ intersecting circles - see the bottom of Figure \ref{tangle}), and 
\item pieces that are not quadratic diagrams and do not have any cutting curves. 
\end{itemize}
These latter pieces are called {\em hyperbolic pieces}.

  The first hyperbolic example is the octahedron (see the two  diagrams on the right in Figure \ref{bye}). 
 The next hyperbolic example is the $3\times 4$ turkshead, with $8$ vertices.  The $3\times n$ turkshead 
 are the diagrams for a standard braid, going around in a circle with $2n$ crossings.  
 
   The $3\times 1$ turkshead is $\{|z|=1\pm \ep\}$ union an arc  in the intersection of the upper half plane with the annulus, connecting $1+\ep$ to $-1+\ep$ and the mirror
 symmetry by complex conjugation, with two vertices $\pm (1+\ep)$. The $3\times n$ turkshead is the inverse image of the $3\times 1$ turkshead by the map $z\mapsto z^n$. 
 This is a lacing pattern for a drum, where a cord zigzags back and forth between a skin on the top and a skin on the bottom and closes in a loop after $2n$ segments. The symmetry group is the dihedral group of order $2n$.

\begin{figure}[h] 
   \centering
   \includegraphics[width=2.5in]{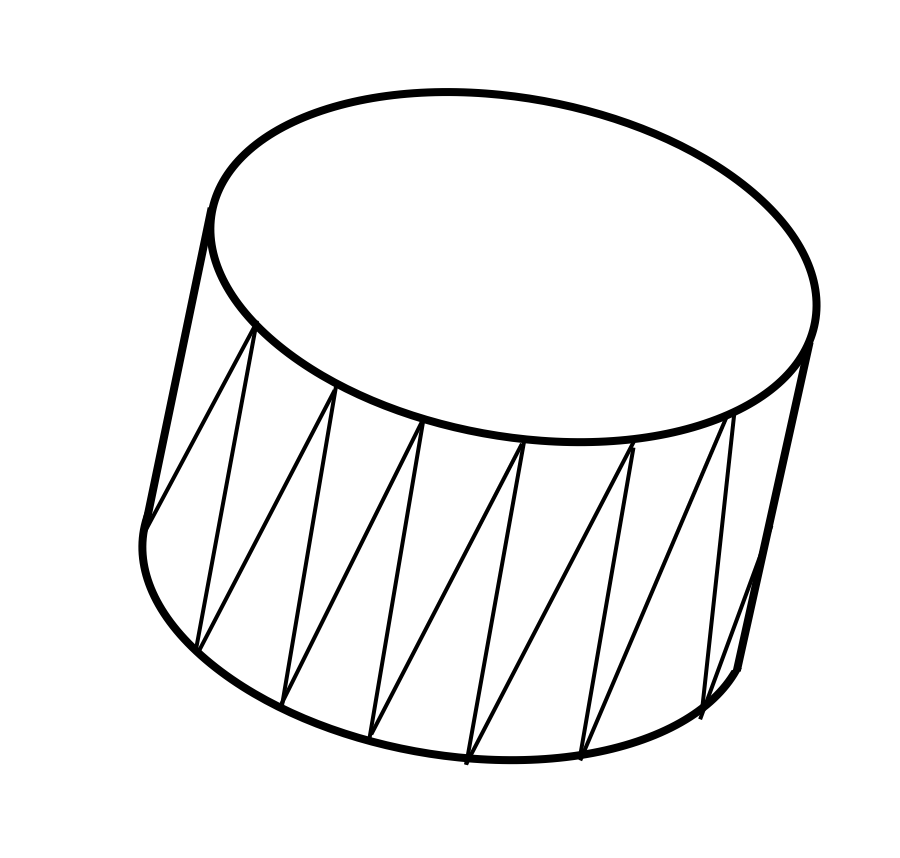} 
   \caption{Turkshead}
   \label{Turkshead}
\end{figure}


  Any 4-valent planar graph is a {\em link projection}: a projection of a link in the 3-ball to the boundary surface in generic positions. There is a nice theory that Martin Bridgman \cite{Br} worked out for the structure of all possible hyperbolic diagrams, although it doesn't take into account the balancing conditions:  they're all generated from $3\times n$ turkshead examples by  choosing an arc cutting across a region between nonadjacent edges but being the first and the third edges of three consecutive edges,  and collapsing the arc to a point, thus creating a new vertex.  If you do this once in a region of each color, the global balance will be retained. 

  \noindent {\bf Problem:} \\
(a) Analyze when the balancing conditions are maintained by these various operations and their inverses.
(This topic has  recently been developed  by J. Tomasini, \cite{To}).\\
(b)  Analyze the set of valid matchings, and the set of valid labelings  under these operations: that is, dots are uniquely partitioned
(or not) to decomposed pieces.

\appendix \section{Hurwitz numbers as seen by Duchi-Poulalhon-Schaeffer}

  Where does the  number  $\dfrac{(2d-2)!}{d!} d^{d-3}$ come from? This is a classical result due to Hurwitz, \cite{H}. There are a number of proofs with different approaches, some applicable to larger classes of Hurwitz numbers. See for example \cite{CT,E1,St,Va}. Here we sketch a  combinatorial approach due to Duchi-Poulalhon-Schaeffer, \cite{DPS}.

  Take a tree (non-embedded) with $d$ white vertices and $d-1$ edges. Label the edges by $1,\cdots, d-1$, in blue.  The number of such edge-labeled trees is known\footnote{Choose
 a root at a  white vertex and label it by $d$. This root choice together with the edge labeling gives a unique labeling of the 
vertices by $1,2,\ldots, d-1,d$ (just label the far end of an edge by the label of the edge). So each root choice 
gives a Cayley tree and the number of such trees is known to be $d^{d-2}$.}  to be $d^{d-3}$.

  For each such tree, put a blue vertex in the middle of
every edge and label it by the original edge label. 
Now label the new set of edges by $1,\ldots, 2d-2$, in red. There are $(2d-2)!$ choices.

  So all together there are $(2d-2)!d^{d-3}$ such bipartite edge-labeled and blue-vertex-labeled trees. 

  Embed  (in a unique way) the tree in the plane such that around each white vertex the increasing blue labels of the incident edges appear clockwise. There are $(2d-2)!d^{d-3}$ such bipartite edge-labeled and blue-vertex-labeled embedded trees. 

  Take a generic branched cover $f:S^2\to S^2$ of degree $d$ with all critical values  on the unit circle, labeled in cyclic order by
$1,\ldots, 2d-2$, in red. Color the point $0$ in blue and the point $\infty$ in white. 
Add in  $2d-2$ disjoint rays (in black) from $0$ to $\infty$ so that each complementary region contains a unique critical value. Call this graph $S$. Now consider the pullback $f^{-1}(S)$ and collapse  each face containing a cocritical point 
to a single edge, obtaining a graph $G$. It has $2d-2$ faces each with a red label of a critical point. This is the dual
graph to underlying $4$-valent graph $f^{-1}(S^1)$. 

  The isomorphism class of $f$ is uniquely determined by this face-labeled graph $G$.

  Now label the blue vertices of $G$ by $1,\cdots,d$ (there are $d!$ such choices). Duchi-Poulalhon-Schaeffer show that
there is a bijection from these face-labeled and blue-vertex-labeled graphs  and the above bipartite edge-labeled and blue-vertex-labeled embedded trees. So the number of isomorphism classes is $\dfrac{(2d-2)!}{d!} d^{d-3}$.

\begin{figure}[h] 
   \centering
   \includegraphics[width=4in]{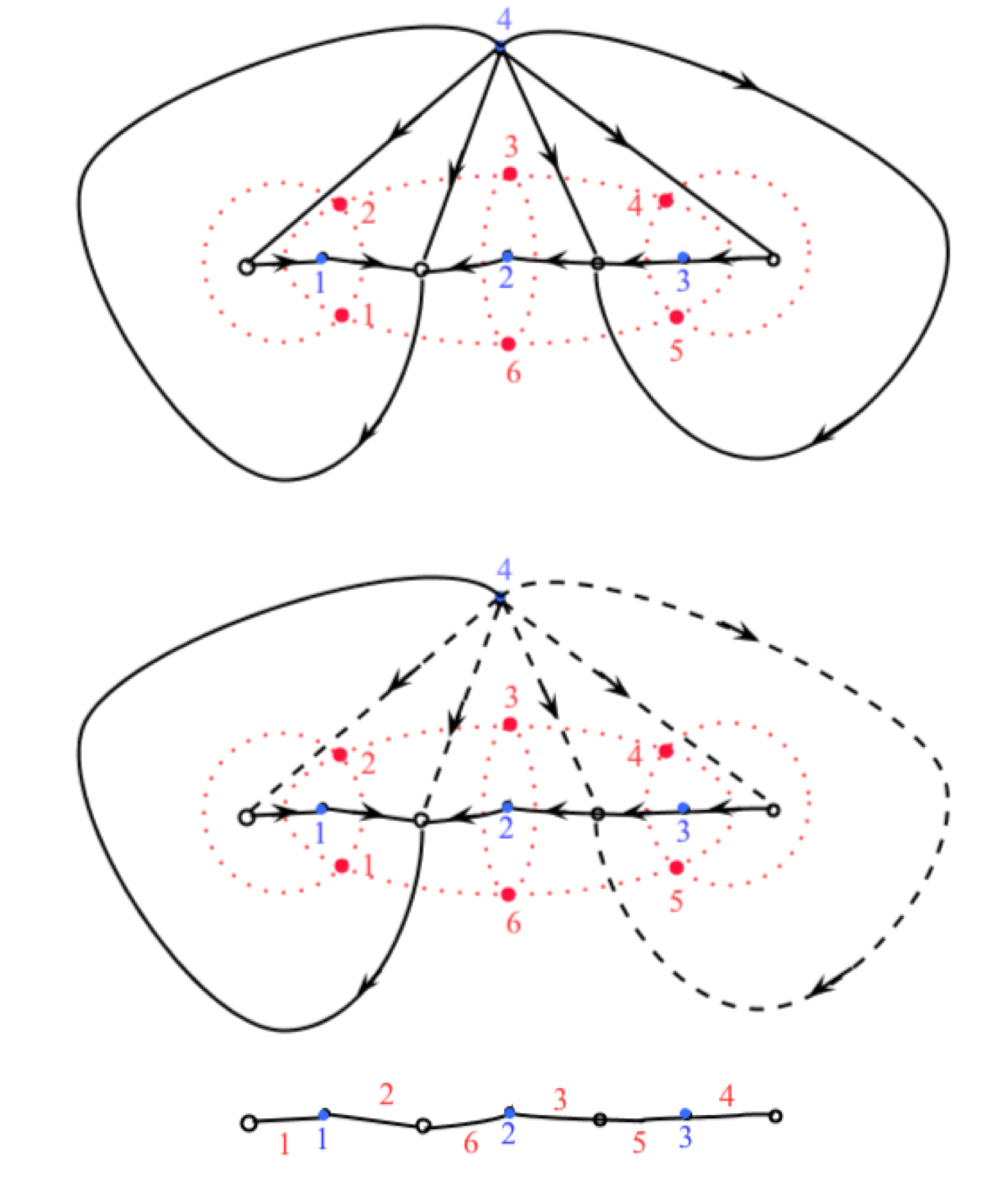} 
   \caption{From $G$ to the spanning tree, and then the to chopped spanning tree}
   \label{span}
\end{figure}

  In fact Duchi, Poulalhon, and Schaeffer take the blue vertex labeled $d$ blue vertex as a root and find a unique spanning tree, as follows:
\begin{itemize}
 \item Embed $G$ in the plane so that the blue vertex with label $d$ is accessible from the unbounded component.  Orient each edge so that the incident face with the greater label is on the left.  Every blue vertex has a unique incoming edge and
every white vertex has a unique outgoing edge. See Figure \ref{span}.

 \item Now  repeatedly reverse orientations of clockwise cycles. This operation does not
change the incoming and outgoing degrees of each vertex. A theorem of Felsner  in \cite{F}  implies that  after finitely many such cycle reversals, we obtain an orientation without clockwise cycles, and this orientation is unique.

 \item Start from the blue vertex labeled $d$  as  root and do a rightmost deep-first search in the sense opposite to the edge orientation, while erasing all other incoming edges to each visited blue vertex (including those incident to the root vertex), to get a canonical spanning tree whose edges are oriented toward the root, as done by Bernardi (\cite{Be}).  See Figure \ref{span}. Chop off the root vertex and its unique incident edge. Label in red each edge of the remaining 
 tree by the label of the white-to-blue right side face. For visual convenience put the label on the same side as the face. Erase the orientation of the tree-edges. See Figure \ref{span}.
\end{itemize}

%

\begin{figure}[h] 
   \centering
   \includegraphics[width=3in]{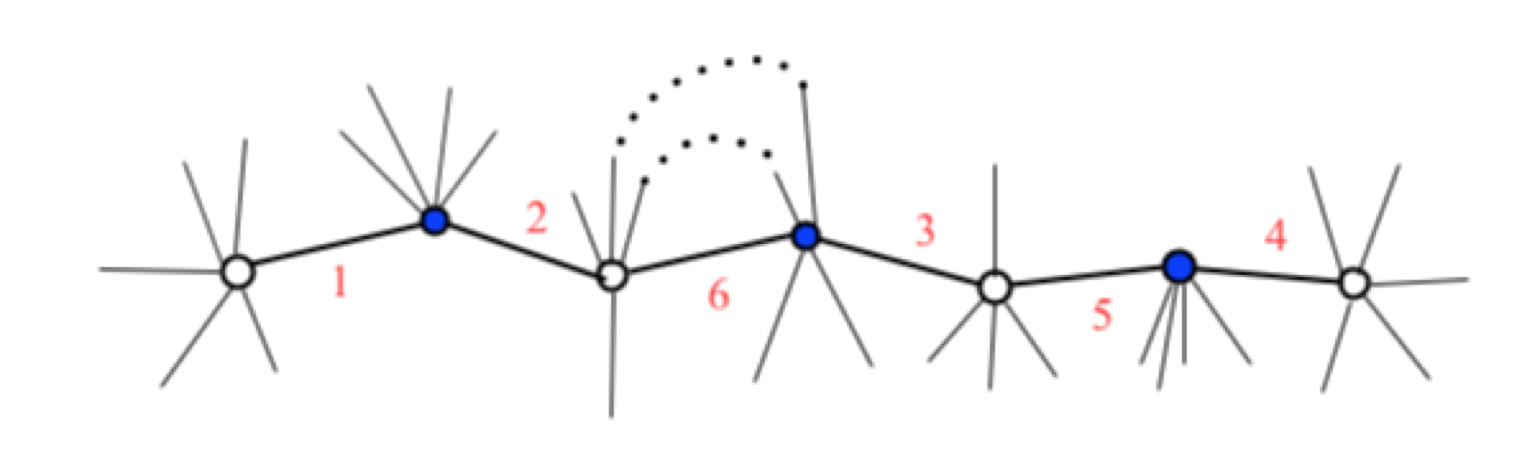} 
   \caption{From a hairy tree and matched hairs (only two steps with dotted arcs are shown) in order to obtain $f^{-1}(\Sigma)$.}
   \label{reverse}
\end{figure}

 The inverse algorithm from an embedded edge-labeled tree to $f^{-1}(\Sigma)$ is also easy to describe; see Figure \ref{reverse}.  Attach some half edges (hairs) to each white and blue vertex to reach the local valence value $2d-2$ and to match the local cyclic order with the labels of the tree edges (counterclockwise for blue and clockwise for white). Now take a blue-incident half edge whose left edge is a full edge, and walk in clockwise order, while matching each visited blue-incident half edge to the closest available white-incident half-edge in the counterclockwise direction. Put an extra blue vertex together with $2d-2$ half edges in the unbounded face and match them  with the remaining $2d-2$ white-incident half-edges.

%

\end{document}